\documentclass[a4paper,11pt]{article}
\usepackage{stmaryrd}
\usepackage{amsfonts}
\usepackage{bbm}
\usepackage{amscd}
\usepackage{mathrsfs}
\usepackage{latexsym,amssymb,amsmath,amscd,amscd,amsthm,amsxtra,xypic}
\usepackage[dvips]{graphicx}
\usepackage[utf8]{inputenc}
\usepackage[T1]{fontenc}
\usepackage{lmodern}
\usepackage{amssymb}
\usepackage[all]{xy}
\usepackage{nicefrac,mathtools,enumitem}
\usepackage{microtype}

\newtheorem{thm}{Theorem}[section]
\newtheorem{lem}[thm]{Lemma}
\newtheorem{cor}[thm]{Corollary}
\newtheorem{pro}[thm]{Proposition}
\newtheorem{ex}[thm]{Example}

\newtheorem{defi}[thm]{Definition}

\newcommand{\be }{\begin{equation}}
\newcommand{\ee }{\end{equation}}

\newcommand{\pf}{\noindent{\bf Proof.}\ }

\def\qed{\hfill ~\vrule height6pt width6pt depth0pt}

\newcommand{\br}[1]{   [ \cdot,    \cdot  ]   }

\newcommand{\Hom}{\mathrm{Hom}}

\newcommand{\Ad}{\mathrm{Ad}}

\newcommand{\End}{\mathrm{End}}

\newcommand{\ISO}{\mathrm{ISO}}

\newcommand{\SO}{\mathrm{SO}}

\begin{document}
\title{\Large\bf
Generalized classical dynamical Yang-Baxter equations and moduli spaces of flat connections on surfaces}

\author{Xiaomeng Xu}

\date{}
\footnotetext{Section of Mathematics, University of Geneva,
2-4 Rue de Li$\rm\grave{e}$vre, c.p. 64, 1211-Gen$\rm\grave{e}$ve
4, Switzerland.\\
E-mail:Xiaomeng.Xu@unige.ch}
%$B$-field transformation,  $E$-dual pair of Lie algebroids}

\maketitle
\begin{abstract}
In this paper, we explain how generalized dynamical r-matrices can be obtained by (quasi-)Poisson reduction. New examples of Poisson structures and Poisson groupoid actions naturally appear in this setting. As an application, we use a generalized dynamical r-matrix induced by the gauge fixing procedure to give a new finite dimensional description of the Atiyah-Bott symplectic structure on the moduli space of flat connections on a surface.
\end{abstract}

\section{Introduction}
The classical Yang-Baxter equation (CYBE) plays a key role in the theory of integrable systems. A geometric interpretation of CYBE was given by Drinfeld and gave rise to the theory of Poisson-Lie groups.
The classical dynamical Yang-Baxter equation (CDYBE) is a differential equation analogue to CYBE and introduced by Felder in \cite{AB} as the consistency condition for the differential Knizhnik-Zamolodchikov-Bernard equations for correlation functions in conformal field theory on tori. It was shown by Etingof and Varchenko \cite{EtingofVarchenko} that dynamical r-matrices correspond to Poisson Lie groupoids (a notion introduced by Weinstein \cite{Alan}) in much the same way as r-matrices correspond to Poisson Lie groups. In the meantime, the classical dynamical Yang-Baxter equation is proven to be closely connected with the theory of Dirac structures and Lie bialgebroids, see \cite{LiuXu} and the references therein. Inspired by the study of Lie bialgebroids, the notion of generalized classical dynamical Yang-Baxter equation was introduced by Liu and Xu \cite{LiuXu} in which the base manifold underlying the CDYBE can be a general Poisson manifold M. Despite its importance, this subject suffered from the lack of examples for a long time.

Since Atiyah and Bott introduced canonical symplectic structures on the moduli spaces of flat connections on Riemann surface in \cite{AtiyahBott}, a lot of attention has been paid to the moduli spaces by mathematicians and physicists due to their rich mathematical structure and their links with a variety of topics. From the physics perspective, a major motivation for their study is their role in Chern-Simons theory.
An independent mathematical motivation for investing moduli spaces of flat connections arises from Poisson geometry.
The Atiyah-Bott symplectic structure on the moduli of flat $G$-connections over oriented surface $\Sigma$ admits several finite dimensional descriptions. The first such description appears in Goldman's study of symplectic structures on character varieties $\Hom(\pi_1(\Sigma),G)/G$, see \cite{Gold}. Another possibility is to obtain the moduli space of flat $G$-connection on a surface $\Sigma_{g,n}$ of genus g with n punctures by (quasi-)Poisson reduction from an enlarged ambient $G^{n+2g}$. In the Fock-Rosly's approach \cite{FockRosly}, the Poisson structure on $G^{n+2g}$ is described using a classical r-matrix. In the theory of Lie group valued moment map \cite{AMM}, the moduli space is obtained by a reduction of a canonical quasi-Poisson tensor on $G^{n+2g}$.

These two lively subjects of dynamical Yang-Baxter equations and moduli spaces of flat connections look very different. However, some recent work indicate the possible connection between them. From the viewpoint of Hamiltonian formalisms, the moduli spaces of flat connections can be viewed as constrained Hamiltonian systems.
Dirac gauge fixing for the moduli space of flat $\ISO(2,1)$-connections on a Riemann surface has been shown to give rise to generalized classical dynamical r-matrices in \cite {Meusburger}.

In this paper, we deepen the connection between these two subjects by giving a systematic investigation of the theory of generalized classical dynamical r-matrices. We explain how generalized dynamical r-matrices can be obtained by (quasi-)Poisson reduction. Further, we show that new examples of Poisson structures and Poisson groupoid actions naturally appear in this setting. After that, we take the canonical quasi-Poisson manifold $G\ast G$ as an example and concretely analyze the dynamical r-matrices arising from the quasi-Poisson reduction of $G\ast G$. We also introduce the notion of gauge transformations for generalized dynamical r-matrices. As an application, we use a dynamical r-matrix induced by the gauge fixing procedure to give a new finite dimensional description of the symplectic structure on the moduli space. We end up with two examples, one of them was previously studied by Meusburger-Sch$\rm\ddot{o}$nfeld in the framework of the $\ISO(2,1)$-Chern-Simons theory of (2+1)-dimensional gravity.

Our paper is structured as follows. In section 2, we recall the definition of generalized classical dynamical Yang-Baxter equations and arising some examples. After that, we give new examples of generalized dynamical r-matrices from (quasi-)Poisson reduction. We show that new Poisson structures and Poisson groupoid actions appear in this setting. In section 3, we use a dynamical r-matrix induced by the gauge fixing procedure to give a new finite dimensional description of the Atiyah-Bott symplectic structure on the moduli space.

\subsection*{Acknowledgements}
\noindent
I give my warmest thanks to my advisor Anton Alekseev for his encouragements as well as inspiring discussions and insightful suggestions. I also want to thank Jianghua Lu for her useful discussions and interest in this paper.
This work is supported by the grant PDFMP2-141756 of the Swiss National Science Foundation.

\section{Generalized classical dynamical Yang-Baxter equations}
Let M be a manifold and $G$ be a Lie group with the Lie algebra $\frak g$. Assume $\{e_i\}_{i=1,..,n}$ is a basis of $\frak g$, then for any tensor $\theta=\sum X_i\otimes e_i\in \Gamma(TM\otimes\frak g)$, smooth map $r:M\rightarrow \frak g\wedge\frak g$ and linear map $\delta:\frak g\rightarrow \frak g\wedge\frak g$, we define the following operation
\begin{eqnarray}
&&\ \ \ \delta \theta=\sum X_i\otimes\delta e_i, \ \ \ \ \ \ \ \ \ \ \ \ \ \ \ \ \ [r,\theta]=\sum X_i\otimes [r,e_i],\label{1}\\
&&[\theta,\theta]=\sum [X_i,X_j]\otimes e_i\wedge e_j, \ \ \ \ \ \ \theta\wedge\theta=\sum X_i\wedge X_j\otimes [e_i,e_j].\label{2}
\end{eqnarray}
Given any $\theta\in \Gamma(TM\otimes\frak g)$, we have the corresponding bundle map $\theta^*:T^*M\rightarrow M\times\frak g$ and a vector bundle map $\theta^{\sharp}:T^*M\rightarrow \frak g$ given by the composition of $\theta^*$ with the projection of $M\times\frak g$ on $\frak g$.
\begin{defi}{\rm{\cite{LiuXu}}}\label{definitionofrmatrix}
For a Poisson manifold $(M,\pi)$ and a Lie algebra $\frak g$, assume that there exists a tensor $\theta\in \Gamma(TM\otimes\frak g)$ such that $\theta^{\sharp}: (T^*M,\pi)\rightarrow \frak g$ is a Lie algebroid morphism. A function $r\in C^{\infty}(M,\wedge^2\frak g)$ is called a dynamical r-matrix coupled with the Poisson manifold $(M,\pi)$ via $\theta$ if
\begin{eqnarray}
\frac{1}{2}[\theta,\theta]=[r,\theta]-\pi^*(dr),
\end{eqnarray}
and the generalized DYBE is satisfied:
\begin{eqnarray}\label{DYBE}
Alt(\theta^* dr)+\frac{1}{2}[r,r]=\Omega,
\end{eqnarray}
where $\Omega\in (\wedge^3\frak g)^{\frak g}$ is an adjoint action invariant element and viewed to be a constant section of $\wedge^3(TM\oplus\frak g)$.
\end{defi}
We call $r$ a triangular dynamical $r$-matrix coupled with $M$ via $\theta$ if the corresponding $\Omega=0$ in \eqref{DYBE}. Throughout this paper, we will also use the triple $(\pi_U,\theta,r)$ to denote a generalized dynamical $r$-matrix.
\begin{ex}\label{classicaldrmatrix}
Let $M=\eta^*$ where $\eta$ is a Lie subalgebra of $\frak g$. Then M is naturally a linear Poisson manifold. Let $\theta^\sharp:T^*\eta^*\rightarrow \frak g$ be the natural projection: $(\xi,v)\rightarrow v$, $\forall~(\xi,v)\in \eta^*\times \eta$. Obviously, $\theta$ is a Lie algebroid morphism. Choose a basis of $\eta$, $\{e_1,...,e_k\}$, and let $\{x_1,...,x_k\}$ be the corresponding coordinate on $\eta^*$. Then $\theta=\sum_i\frac{\partial}{\partial x_i}\otimes e_i$. It is obvious that $[\theta,\theta]=0$. Therefore, the condition $\frac{1}{2}[\theta,\theta]=[r,\theta]-\pi^*(dr)$ takes the form $[r,\theta]=\pi^*(dr)$, which says that the map $r:\eta^*\rightarrow \frak g\wedge \frak g$ is H-equivariant, where $H$ is the connected Lie group with Lie algebra $\eta$.
In this case, equation $Alt(\theta^*dr)+\frac{1}{2}[r,r]=\Omega\in (\wedge^3\frak g)^{\frak g}$ turns to be the classical dynamical Yang-Baxter equation and a solution $r$ is a classical dynamical r-matrix \cite{AB}.
\end{ex}
Similarly, we can introduce a notion of the generalized Poisson-Lie dynamical Yang-Baxter equation.
\begin{defi}\label{PoissonLie}
For a Poisson manifold $(M,\pi)$ and a Lie bialgebra $(\frak g,\delta)$, assume that there exists a tensor $\theta\in \Gamma(TM\otimes\frak g)$ such that $\theta^{\sharp}: (T^*M,\pi)\rightarrow \frak g$ is a Lie algebroid morphism. A function $r\in C^{\infty}(M,\wedge^2\frak g)$ is called a Poisson-Lie dynamical r-matrix coupled with the Poisson manifold $(M,\pi)$ via $\theta$ if
\begin{itemize}
\item[(a)] $\delta\theta+\frac{1}{2}[\theta,\theta]=[r,\theta]-\pi^{\sharp}(dr),$

\item[(b)] the generalized Poisson-Lie DYBE is satisfied:
\begin{eqnarray}\label{PDYBE}
Alt(\theta^{*}dr)+\frac{1}{2}[r,r]+\delta r=\Omega\in (\wedge^3\frak g)^{\frak g}.
\end{eqnarray}
\end{itemize}
\end{defi}
\begin{ex}
Let $(\frak g,\frak g^*)$ be a Lie bialgebra, $G$ and $G^*$ be the simply connect Lie group corresponding to $\frak g$ and $\frak g^*$ respectively. There is a bivector field $\pi$ such that $(G^*,\pi)$ is a Poisson Lie group. A natural section $\theta$ of $TG^*\otimes \frak g$ is induced by the isomorphism $T^*G^*\rightarrow G\times\frak g$. It is easy to check that $\theta^\sharp: T^*G^*\rightarrow \frak g$ is a Lie algebroid morphism. A direct calculation shows that equation $(a)$ and $(b)$ in Definition \ref{PoissonLie} can be written as

$$dress_a^L(r)+[a\otimes 1+1\otimes a,r]=0,$$
$$[r,r]+Alt(d^L r)+Alt((\delta\otimes id)(r))=\Omega,$$
for $r:G^*\rightarrow \frak g\wedge\frak g$ and any $a\in \frak g$, where $dress_a^L$ denotes the left dressing vector field generalized by $a$ and $d^Lr(g)=e_i\otimes \frac{d}{dt}|_{t=0}r(exp(-te_i)\cdot g)$ for any $g\in G^*$ and a basis $\{e_i\}$ of $\frak g$. Thus a map $r:G^*\rightarrow \frak g\wedge\frak g$ is a generalized dynamical r-matrix coupled with $G^*$ via $\theta$ if r is a Poisson-Lie dynamical r-matrix \cite{EnriquezEtingof}.
\end{ex}

\subsection{Generalized classical dynamical r-matrices from (quasi-)Poisson reduction}
In this subsection, we show that generalized classical dynamical r-matrices naturally appear in the theory of (quasi-)Poisson reduction. First, let us recall the definition of quasi-Poisson $G$-manifolds.

We assume $(\frak g,\langle \cdot,\cdot\rangle )$ is a quadratic Lie algebra, $\phi$ is the Cartan 3-tensor. In terms of an orthogonal basis $\{e_a\}$ of $\frak g$, $\phi\in\wedge^3\frak g$ is given by
\begin{eqnarray}
\phi=\frac{1}{12}f_{abc}e_a\wedge e_b\wedge e_c,
\end{eqnarray}
where $f_{abc}=\langle e_a,[e_b,e_c]\rangle$ are the structure constants of $\frak g$.
For any $G$-manifold $M$, the Lie algebra homomorphism $\rho:\frak g\rightarrow TM$ can be extended to an equivariant map,
$$
\rho:\wedge^{\bullet}\frak g\rightarrow \wedge^{\bullet}TM,
$$
preserving wedge products and Schouten brackets.

\begin{defi}{\rm \cite{Anton}}
A quasi-Poisson manifold is a $G$-manifold $M$, equipped with an invariant bivector field $\pi\in \Gamma(\wedge^2TM)$ such that
\begin{eqnarray}
[\pi,\pi]=\rho(\phi).
\end{eqnarray}
\end{defi}
\begin{ex}\label{quasiG}
Let $G$ be a Lie group and $\{e_a\}_{a\in I}$ be an orthogonal basis of its Lie algebra $\frak g$. Define a bivector on $G$ by $\pi_G=\frac{1}{2}\sum_{a\in I} R_a\wedge L_a$, where $R_a$ and $L_a$ are right and left invariant vectors generated by $e_a$. Then $(G,\pi_G)$ is a quasi-Poisson $G$-manifold, where $G$ acts on itself by conjugation.
\end{ex}

Generally, the $G$-invariant functions on a quasi-Poisson manifold $(M,\pi)$ is a Poisson algebra under the binary bracket induced by $\pi$. Thus we can get a Poisson algebra on $C^{\infty}(M)^G$.

Let $M$ be a $G$-manifold with $G$ acting locally freely and $\rho:\frak g\rightarrow TM$ be the corresponding infinitesimal action. We use the same symbol $\rho$ to denote the following natural extension map:
\begin{eqnarray}
\rho:\wedge^{\bullet}(TM\oplus\frak g)\rightarrow \wedge^{\bullet}TM.
\end{eqnarray}
Throughout this paper, we denote the anti-symmetrization of any section $A\in\Gamma(\wedge^{\bullet}(TM)\otimes\wedge^{\bullet}\frak g)$ by $\hat{A}\in\Gamma(\wedge^{\bullet}(TM\oplus\frak g))$. Thus, if $\theta=f^{ia}(x)\frac{\partial}{\partial x_i}\otimes e_a\in\Gamma(TM\otimes\frak g)$ in local coordinates $\{x_i\}$ of $U$ and a basis $\{e_a\}$ of $\frak g$, we have $\hat{\theta}=f^{ia}(x)\frac{\partial}{\partial x_i}\wedge e_a$.
\begin{thm}\label{ReduceGDYBE}
Let $U\subset M$ be a cross-section of the $G$ action and $\pi_M$ be a bivector field on $M$. Then there exists a unique triple $(\pi_U, \theta,r)$, where $\pi_U\in\Gamma(\wedge^2 TU)$, $\theta\in\Gamma(TU\otimes\frak g)$ and $r:U\rightarrow \frak g\wedge\frak g$ such that
\begin{eqnarray}\label{restrictionpi}
\pi_M|_U=\pi_U-\rho(\hat{\theta})|_U+\rho(r)|_U.\label{-}
\end{eqnarray}
Moreover,
\begin{itemize}
\item[(a)] if $\pi_M$ is a $G$-invariant Poisson tensor on $M$, then $(U,\pi_U)$ is a Poisson manifold and $r$ is a triangular dynamical r-matrix coupled with $U$ via $\theta$.

\item[(b)] if $\pi_M$ is a quasi-Poisson tensor on $M$, then $(U,\pi_U)$ is a Poisson manifold and $r$ is a dynamical r-matrix coupled with $U$ via $\theta$ with respect to the Cartan 3-tensor, i.e., $\Omega=-\frac{1}{2}\phi$ in the generalized CDYBE.
\end{itemize}
\end{thm}
\pf Because $G$ acts locally freely and $U$ is a cross-section, for any $x \in U$ there is a canonical splitting $T_xM= T_xU \oplus \rho_x(\frak g)$ of the sequence $$0\rightarrow\frak g \rightarrow T_xM \rightarrow T_xU\rightarrow 0.$$ Thus, there exists unique $\pi_U\in \Gamma(\wedge^2TU)$, $\theta\in\Gamma(TU\otimes \frak g)$ and $r:U\rightarrow \frak g\wedge\frak g$ such that $$\pi_M|_U=\pi_U-\rho(\hat{\theta})|_U+\rho(r)|_U,$$ where $\pi_U$ is tangent to $U$.

If $\pi_M$ is a $G$-invariant Poisson tensor or quasi-Poisson tensor, it induces a Poisson bracket $\{\cdot,\cdot\}$ on $U$. From the expression \eqref{restrictionpi} and the fact $\rho(e)f'=0$ for any $e\in\frak g$ and $f'\in C^{\infty}(M)^G$, we have that $\{f,g\}=\pi_U(df,dg)$ for any $f,g\in C^{\infty}(U)$. This is to say $(U,\pi_U)$ is a Poisson manifold. To complete the proof, we have to check that the required compatibility condition and the generalized CDYBE are satisfied by the triple $(\pi_U,\theta,r)$. This proof is similar to the ones of Theorem \ref{PoissonDYBE} and Theorem \ref{ReducedPoisson}.
\qed
\\

Theorem \ref{ReduceGDYBE} suggests the following construction which generalizes the construction for ordinary classical dynamical r-matrices in \cite{PingXu}. Given a manifold $M$, $M\times G$ carries natural right and left $G$-actions defined by $(x,p)\cdot g=(x,pg)$ and $g\cdot (x,p)=(x,gp)$ respectively, for all $x\in M$, $p,g\in G$. Then we have:
\begin{thm}\label{PoissonDYBE}
Let $(M,\pi_M)$ be a Poisson manifold and $\theta\in\Gamma(TM\otimes\frak g)$. Then any smooth function $r:M\rightarrow \frak g\wedge \frak g$ induces a right $G$-invariant bivector $\pi_r$ on $M\times G$ which is given by
\begin{eqnarray}
\pi_r=\pi_M+\rho^L(\hat{\theta})+\rho^L(r),\label{+}
\end{eqnarray}
and

$(a)$ $\pi_r$ is a Poisson tensor iff $r$ is a triangular generalized dynamical $r$-matrix.

$(b)$ $\pi_r$ is a quasi-Poisson tensor iff $r$ is a generalized dynamical $r$-matrix with respect to the Cartan 3-tensor.
\end{thm}
\pf
Note that the vector field on $M\times G$ has a natural bigrading: elements in $TM$ have degree $(1,0)$ while elements in $TG$ have degree $(0,1)$. It is simple to see that $[\pi_M,\pi_M]$ is of degree $(3,0)$, $[\pi_M,\rho^L(\hat{\theta})]$ is of degree $(2,1)$, $[\pi_M,\rho^L(r)]$ is of degree $(1,2)$ and $[\rho^L(r),\rho^L(r)]$ is of degree $(0,3)$. On the other hand, $[\rho^L(\hat{\theta}),\rho^L(r)]$ consists of elements of degree $(1,2)$ and of $(0,3)$ and $[\rho^L(\hat{\theta}),\rho^L(\hat{\theta})]$ consists of elements of degree $(2,1)$ and of $(1,2)$. For any $S\in \wedge^3(TM\oplus TG)$, let $S=\sum_{0\le i,j\le 3}S^{(i,j)}$ be its decomposition with respect to this bigrading. The following equations can be verified by a direct computation:
\begin{eqnarray}
&&[\rho^L(\hat{\theta}),\rho^L(\hat{\theta})]^{(1,2)}=\rho^L(\widehat{[\theta,\theta]}), \ \ \ \ \ \ \ \ \ \ [\rho^L(\hat{\theta}),\rho^L(\hat{\theta})]^{(2,1)}=2\rho^L(\widehat{\theta\wedge\theta}),\\
&&[\rho^L(\hat{\theta}),\rho^L(r)]^{(0,3)}=\rho^L(Alt(\theta^*dr)), \ \ \ \ [\pi_M,\rho^L(r)]=\rho^L(\pi_M^{\sharp}(dr)),\\
&&[\rho^L(\hat{\theta}),\rho^L(r)]^{(1,2)}=-\rho^L(\widehat{[r,\theta]}), \ \ \ \ \ \ \ \ \ [\pi_M,\rho^L(\hat{\theta})]=\rho^L(d_{\pi_M}\theta)
\end{eqnarray}
where the operations $[\theta,\theta]$, $\theta\wedge\theta$ and $[r,\theta]$ for $\theta\in\Gamma(TM\otimes\frak g)$ and $r\in C^{\infty}(M,\wedge^2\frak g)$ are defined as \eqref{1} and \eqref{2}. Eventually, by using the facts $[\rho^L(r),\rho^L(r)]=\rho^L([r,r])$ and $[\pi_M,\pi_M]=0$ we have
\begin{eqnarray*}
[\pi_r,\pi_r]&=&[\pi_M+\rho^L(\hat{\theta})+\rho^L(r),\pi_M+\rho^L(\hat{\theta})+\rho^L(r)]\nonumber\\
&=&[\rho^L(\hat{\theta}),\rho^L(\hat{\theta})]+2[\pi_M,\rho^L(\hat{\theta})]+2[\rho^L(\hat{\theta}),\rho^L(r)]+2[\pi_M,\rho^L(r)]+[\rho^L(r),\rho^L(r)]\nonumber\\&=&2\rho^L(Alt(\theta^*dr)+\frac{1}{2}[r,r])+2\rho^L(\widehat{[\theta,\theta]}-\widehat{[r,\theta]}+\widehat{\pi^*(dr)})+2\rho^L(\widehat{\theta\wedge\theta}+d_{\pi_M}\theta).
\end{eqnarray*}
Note that the map $\theta^\sharp:T^*M\rightarrow \frak g$ is a Lie algebroid morphism is equivalent to $$\widehat{\theta\wedge\theta}+d_{\pi_M}\theta=0.$$ Therefore we have that $[\pi_r,\pi_r]=0$ iff $r$ is a triangular generalized dynamical $r$-matrix and $[\pi_r,\pi_r]=\rho^R(\phi)$ iff $r$ is a generalized dynamical $r$-matrix with respect to $\Omega=-\frac{1}{2}\phi$.
\qed
\\

Similarly, we have the following theorem.
\begin{thm}\label{adjointPoisson}
Let $(N,\pi_N)$ be a quasi-Poisson G-space and $\rho:\frak g\rightarrow TN$ be the infinitesimal action. Then for any generalized dynamical r-matrix coupled with $(M,\pi_M)$ via $\theta$ with respect to $\Omega=-\frac{1}{2}\phi$,
\begin{eqnarray}
\pi:=\pi_M+\pi_N+\rho(\hat{\theta})+\rho(r)
\end{eqnarray}
is a Poisson tensor on $M\times N$.
\end{thm}
\pf We need to prove $[\pi,\pi]=0$. Note that $[\rho(r),\pi_N]=[\rho(\theta),\pi_N]=0$ because of the invariance of $\pi_N$. Thus we have
\begin{eqnarray*}
[\pi,\pi]&=&[\pi_M+\pi_N+\rho(\hat{\theta})+\rho(r),\pi_M+\pi_N+\rho(\hat{\theta})+\rho(r)]\nonumber\\&=&[\pi_M,\pi_M]+[\pi_N,\pi_N]+\rho([r,r])+2[\rho(\hat{\theta}),\rho(r)]+2[\pi_M,\rho(\hat{\theta})]+2[\pi_M,r]\nonumber\\&=&[\pi_N,\pi_N]+2\rho(Alt(\theta^* dr)+[r,r])\\
&=&\rho(\phi)-\rho(\phi)=0.
\end{eqnarray*}
This finishes the proof.
\qed
\\

Now we discuss the relation between the generalized dynamical r-matrix and the reduction of the fusion product of two quasi-Poisson manifolds. Let $M$, $N$ be two G-manifolds and $\rho_M$, $\rho_N$ be the corresponding infinitesimal $G$ action. We define a bivector field on $M\times N$ by $$\Phi=\sum_{a\in I}\rho_M(e_a)\wedge \rho_N(e_a),$$ where $\{e_a\}_{a\in I}$ is an orthogonal basis of $\frak g$.
\begin{pro}{\rm \cite{Anton}}
If $(M,\pi_M)$ and $(N,\pi_N)$ are two quasi-Poisson G-manifolds, then $\pi_M+\pi_N-\Phi$ gives a quasi-Poisson structure on $M\times N$ for the diagonal $G$-action. This quasi-Poisson manifold, denoted by $M\ast N$, is called the fusion product of $M$ and $N$.
\end{pro}
\begin{ex}\label{quasiG2}
Let $(G,\pi_G)$ be the quasi-Poisson $G$-manifold given in Example \ref{quasiG}. By using the fusion product between $G$ with itself, we get a quasi-Poisson manifold $G\ast G$. Let $R^i_a$ and $L^i_a$ denote the right and left invariant vector fields on $i$-th copy of $G\times G$ generated by $e_a$, then the quasi-Poisson tensor takes the form:
\begin{eqnarray}
\pi_{G^2}=\frac{1}{2}\sum_{a} (R^1_a\wedge L^1_a+R^2_a\wedge L^2_a+(L^1_a-R^1_a)\wedge (L^2_a-R^2_a)).
\end{eqnarray}
\end{ex}

Let $(M,\pi_M)$ and $(N,\pi_N)$ be two quasi-Poisson $G$-manifolds. We assume $G$ acts locally freely on $M$. Let $U\subset M$ be any cross-section of the $G$ action on $M$. By Theorem \ref{ReduceGDYBE}, associated to $U$, there is a triple $(\pi_U, \theta,r)$ such that $r$ is a generalized dynamical r-matrix coupled with $(U,\pi_U)$ via $\theta$. On the other hand, $U\times N$ is a cross-section of the diagonal $G$ action on $M\times N$. So it inherits a Poisson structure $\pi_{\rm red}$ by reduction from the quasi-Poisson structure on $M\ast N$.
\begin{thm}\label{ReducedPoisson}
The Poisson tensor $\pi_{red}$ on $U\times N$ takes the form
\begin{eqnarray}
\pi_{{\rm red}}=\pi_U+\rho_N(\hat{\theta})+\rho_N(r)+\pi_N.
\end{eqnarray}
\end{thm}
Before giving a proof, we show the following lemma which says that in the reduction level, fusion product and direct product are same.

\begin{lem}\label{reductionlevel}
Let $M$ and $N$ be two quasi-Poisson G-manifolds. Then for any diagonal $G$-invariant functions $f, g\in C^{\infty}(M\times N)^G$, $\Phi(df,dg)=0$. Moreover, $\pi_M+\pi_N$ induces a Poisson algebra structure on $C^{\infty}(M\times N)^G$ which is same as the Poisson algebra on $C^{\infty}(M\ast N)^G$.
\end{lem}
\pf Note that $\Phi=\sum \rho_M(e_a)\wedge \rho_N(e_a)$, and $(\rho_M(e_a)+\rho_N(e_a))f=0$ for $f$ $G$-invariant function on $M\times N$. So $\Phi=-\sum \rho_N(e_a)\wedge \rho_N(e_a)=0$ when restricts to G-invariant functions. \qed

{\bf Proof of Theorem \ref{ReducedPoisson}} For any $f,g\in C^{\infty}(U\times N)$, let $f',g'\in C^{\infty}(M\times N)^G$ be the diagonal $G$-invariant extension of $f$, $g$ respectively. Then by Lemma \ref{reductionlevel}, $\pi_{{\rm red}}(df,dg)=(\pi_M+\pi_N)(df',dg')|_{U\times N}$. Following Theorem \ref{ReduceGDYBE}, we have $$\pi_M|_{U}=\pi_U-\rho_M(\hat{\theta})|_U+\rho_M(r)|_U.$$ Together with the fact $(\rho_M(e_a)+\rho_N(e_a))F=0$ for any $F\in C^{\infty}(M\times N)^G$, we get $$(\pi_M+\pi_N)(df',dg')|_{U\times N}=(\pi_U+\rho_N(\hat{\theta})+\rho_N(r)+\pi_N)(df',dg')|_{U\times N}.$$ Note that $\pi_U+\rho_N(\hat{\theta})+\rho_N(r)+\pi_N$ is tangent to $U\times N$ and $f'|_U=f$, $g'|_U=g$, therefore, $$(\pi_U+\rho_N(\hat{\theta})+\rho_N(r)+\pi_N)(df',dg')|_{U\times N}=(\pi_U+\rho_N(\hat{\theta})+\rho_N(r)+\pi_N)(df,dg).$$ This is to say $\pi_{{\rm red}}=\pi_U+\rho_N(\hat{\theta})+\rho_N(r)+\pi_N$. \qed

\subsection{Generalized dynamical r-matrices from the reduction of $(G\times G,\pi_{G^2})$}

Let $(G\times G,\pi_{G^2})$ be the quasi-Poisson manifold given in Example \ref{quasiG2}. For a conjugacy class $\mathcal{C}$ in G, identify the tangent space $T_g\mathcal{C}$ at $g$ with $\frak g_g^{\bot}$, where $\frak g_g$ is the Lie algebra of the stabilizer of $\frak g$. The operator $Ad_g-1$ in invertible on $\frak g_g^{\bot}$. Thus we get a linear operator
$$
\frac{\Ad_g+1}{\Ad_g-1}|\frak g_g^{\bot}:=(\frac{\Ad_g+1}{\Ad_g-1})\circ {\rm Pr}_{\frak g_g^\bot}:\frak g \rightarrow \frak g,
$$
where $\Pr_{\frak g_g^\bot}$ is the projection of $\frak g$ on $\frak g_g^\bot$.
\begin{pro}[Proposition 3.4, \cite{Anton}]\label{conjugacyPoisson}
$\sum_{a\in I} R_a\wedge L_a=\sum_{a,b\in I}\frac{1}{2}(\frac{\Ad_g+1}{\Ad_g-1}|\frak g_g^{\bot})_{ab}X_a\wedge X_b$ as bivector fields on $G$, where $\{e_a\}_{a\in I}$ is a basis of $\frak g$ and $X_a=L_a-R_a$ for any $e_a\in\frak g$.
\end{pro}
As a result, for any conjugacy classes $\mathcal{C}_1$, $\mathcal{C}_2$ in $G$, ${\pi_{G^2}|}_{{\mathcal{C}_1\times \mathcal{C}_2}}$ is tangent to $\mathcal{C}_1\times \mathcal{C}_2$, i.e., $(\mathcal{C}_1\times\mathcal{C}_2,{\pi_{G^2}|}_{\mathcal{C}_1\times\mathcal{C}_2})$ is a quasi-Poisson manifold.

Now let $T\subset G$ be a maximal torus, $p\in \mathcal{C}_1\cap T$ and $G_p$ its isotropic group under the conjugation action of G, $\frak g_p$ the Lie algebra of $G_p$. For any open set in $\mathcal{C}_2$ where the conjugation $G_p$ action is locally free, we can choose a cross-section $U$ of the $G_p$ action. Thus, we have $\{p\}\times U$ is a cross-section of the simultaneous conjugation G action on $\mathcal{C}_1\times \mathcal{C}_2$. Following Theorem \ref{ReduceGDYBE}, associated to $\{p\}\times U$, there is a dynamical r-matrix $(\pi_{p\times U},\theta,r)$. To write down it explicitly, we introduce a function $H\in C^{\infty}(U, \End(\frak g))$ in the following way. For any point $x\in U$, let $\frak g'_x$ be the subspace of $\frak g$ defined by $$\frak g'_x=\{u\in \frak g|\frac{d}{dt}|_{t=0}exp(tu)\cdot x\cdot exp(-tu)\in T_xU\}.$$ Because $U$ is a cross-section of the conjugation $G_p$ action, we have a set of direct product decompositions of $\frak g$ as $\frak g=\frak g_p\oplus \frak g'_x$ parameterized by the coordinates on $U$. Then we define $H_x\in\End(\frak g)$ to be the projection of $\frak g$ on $\frak g'_x$. For any $e_a\in\frak g$, we have $H(e_a)\in C^{\infty}(U,\frak g)$ and a vector field $H(X_a)\in\Gamma(TU)$ defined by
\begin{eqnarray}
&&H(X_a)|_x:={\frac{d}{dt}|}_{t=0}exp(tH_x(e_a))\cdot x\cdot exp(-tH_x(e_a)), \ \forall x\in U.
\end{eqnarray}

\begin{thm}\label{conjugacyclassrmatrix}
The triple $(\pi_{p\times U},\theta,r)$ from the reduction of the quasi-Poisson tensor $\pi_{G^2}$ on $\mathcal{C}_1\times\mathcal{C}_2$ is given by
\begin{eqnarray}
&&\pi_{p\times U}=\frac{1}{2}\sum_{a,b}((\frac{\Ad_p+1}{\Ad_p-1}|_{\frak g_p^\bot})_{ab}+\frac{\Ad_x+1}{\Ad_x-1}|\frak g_x^{\bot})_{ab})H(X_a)\wedge H(X_b),\\
&&\theta=\frac{1}{2}\sum_{a} H(X_a)\otimes e_a+\frac{1}{2}(\frac{\Ad_p+1}{\Ad_p-1}|_{\frak g_p^\bot})_{ab}H(X_a)\otimes H(e_b)+\nonumber\\&& \ \ \ \ \ \frac{1}{2}\sum_{a,b}(\frac{\Ad_x+1}{\Ad_x-1}|\frak g_x^{\bot})_{ab}H(X_a)\otimes (H(e_b)-e_b),\\
&& r=\frac{1}{2}\sum_{a} e_a\wedge H(e_a)+\frac{1}{2}(\frac{\Ad_p+1}{\Ad_p-1}|_{\frak g_p^\bot})_{ab}H(e_a)\wedge H(e_b)+\nonumber\\&& \ \ \ \ \ \frac{1}{2}\sum_{a,b}(\frac{\Ad_x+1}{\Ad_x-1}|\frak g_x^{\bot})_{ab}(e_a-H(e_a))\wedge (e_b-H(e_b)).
\end{eqnarray}

\end{thm}
\pf Note that for the cross-section $p\times U$ of the $G_p$ action on $\mathcal{C}_1\times\mathcal{C}_2$, the triple $(\pi_{p\times U},\theta,r)$ corresponds to the decomposition of $\pi_{G^2}|_{p\times U}$ with respect to $TU$ and the complement $\rho(\frak g)|_U$ generated by the diagonal $G$-action, where $\rho(e_a)=X^1_a+X^2_a$ for any $e_a\in\frak g$. On the other hand, from Example \ref{quasiG2} and Proposition \ref{conjugacyPoisson}, we have at any point $(p,x)\in \{p\}\times U$,
\begin{eqnarray}\label{piG2}
\pi_{G^2}=\frac{1}{4}\sum_{a,b}(\frac{\Ad_p+1}{\Ad_p-1}|\frak g_p^{\bot})_{ab}X^1_a\wedge X^1_b+\frac{1}{4}\sum_{a,b}(\frac{\Ad_x+1}{\Ad_x-1}|\frak g_x^{\bot})_{ab}X^2_a\wedge X^2_b+\sum_{a}X_a^1\wedge X_a^2.
\end{eqnarray}
With the help of this expression, we just need to compute the decomposition of $X^i_a\in \Gamma(T(\mathcal{C}_1\times\mathcal{C}_2)|_U)$ along the two directions $TU$ and $\rho(\frak g)|_U$. After a direct computation, we get the following decompositions:
\begin{eqnarray}
&&X_a^1|_{p\times U}=-H(X_a)+\rho(H(e_a))|_{p\times U},\label{decomposition1}\\
&&X_a^2|_{p\times U}=H(X_a)+\rho(e_a-H(e_a))|_{p\times U},\label{decomposition2}
\end{eqnarray}
where $\rho(H(e_a))\in \Gamma(T(\mathcal{C}_1\times\mathcal{C}_2))$ is given by $\rho(H(e_a))|_{y,x}=(X^1_{H_x(e_a)}+X^2_{H_x(e_a)})|_{y,x}$ for all $(y,x)\in \mathcal{C}_1\times\mathcal{C}_2$.
To be precise, by the definition of $H$, $e_a-H(e_a)\in C^{\infty}(U,\frak g_p)$ where $\frak g_p$ is the Lie algebra of the isotropic group $G_p$, so we get $X^1_a=X^1_{H(e_a)}$ when restricts to $p\times U$. Thus we have $$-H(X_a)+\rho(H(e_a))|_{p\times U}=-H(X_a)+X^1_{H(e_a)}|_{p\times U}+H(X_a)=X_a^1|_{p\times U}.$$ A similar calculation gives the equation \eqref{decomposition2}.

In the end, we get the expression of the triple $(\pi_{p\times U},\theta,r)$ by plugging \eqref{decomposition1} and \eqref{decomposition2} in the expression of $\pi_{G^2}|_{p\times U}$.
\qed
\\

We refer to the generalized dynamical r-matrices associated to two conjugacy classes in $G$ as moduli space generalized dynamical r-matrices.

Let us take $G=\rm{SU(2)}$ for a concrete example. Let
$$e_1=
{\footnotesize\left(
  \begin{array}{cc}
    0 & i  \\
    i & 0
  \end{array}
\right)}
, \quad
e_2=
{\footnotesize\left(
  \begin{array}{cc}
    i & 0  \\
    0 & -i
  \end{array}
\right)}
,\quad
e_3=
{\footnotesize\left(
  \begin{array}{cc}
    0 & 1 \\
    -1 & 0
  \end{array}
\right)}
$$
be a basis of $su(2)$ and
$\mathcal{C}\subset \rm{SU(2)}$ the conjugacy class through p=
$
{\footnotesize\left(
  \begin{array}{cc}
    i & 0  \\
    0 & -i
  \end{array}
\right)}
$.
Then $\mathcal{C}$ can be identified with the sphere $S^2=\{(x,y,z)\in\mathbb R^3:x^2+y^2+z^2=1\}$ and an element in $\mathcal{C}$ takes the form
$
{\footnotesize\left(
  \begin{array}{cc}
    ix & y+iz \\
    -y+iz & -ix
  \end{array}
\right)}
$.
The diagonal matrix
$
{\footnotesize\left(
  \begin{array}{cc}
    e^{i\beta} & 0 \\
    0 & e^{-i\beta}
  \end{array}
\right)}
$
acts on $\mathcal{C}=S^2\in \mathbb R^3$ by rotation with respect to $x$-axis, i.e., $e^{i\beta}\circ (x,y,z)=(x,e^{2i\beta}y,e^{2i\beta}z)$. For this $S^1$ action, we choose a simple cross-section $$U:=\{(x,y,z)~|~-1< x< 1,~y=0,~z>0\}\subset \mathcal{C}=S^2\in\mathbb R^3.$$ $U$ can be parameterized by $\alpha$ where $x=\sin\alpha$ and $z=\cos\alpha$. Recall that at any point $\alpha\in U$, we have the subspace of $\frak g$ defined by $$\frak g'_{\alpha}=\{e\in\frak g|\frac{d}{dt}_{|_{t=0}}exp(te)\cdot \alpha \cdot exp(-te)\in T_{\alpha}U\}.$$ A direct calculation gives
\begin{pro}
In terms of $e_1$, $e_2$, $e_3$, $\frak g'_{\alpha}={\rm Span}\{e_1+\tan\alpha e_2,e_3\}$.
\end{pro}

It follows that the function $H\in C^{\infty}(U,\End(\frak g))$ is given by
$$
H(e_1)=0, \ H(e_2)=\cot \alpha e_1+e_2, \ H(e_3)=e_3.
$$
The corresponding vector fields on $U$ generated by adjoint action are given by
$$
H(X_1)=0, \ H(X_2)=0, \ H(X_3)=\frac{\partial}{\partial\alpha}.
$$
Another straightforward computation shows that $\frak g_{\alpha}^{\bot}={\rm Span}\{\tan \alpha e_1-e_2,e_3\}$, where $\frak g_{\alpha}$ is the Lie subalgebra of the stabilizer of $\frak g$ at $\alpha\in U$. However, we have $$(\Ad_{\alpha}+1)e_3=(\Ad_{\alpha}+1)(\tan \alpha e_1-e_2)=0.$$ It is to say that $(\frac{\Ad_g+1}{\Ad_g-1}|\frak g_{\alpha}^{\bot}):\frak g_{\alpha}^{\bot}\rightarrow \frak g_{\alpha}^{\bot}$ is a zero transformation.

Eventually, by Theorem \ref{conjugacyclassrmatrix}, the dynamical r-matrix associated to the local section ${p}\times U$ of $(\mathcal{C}\times\mathcal{C})/G$ is given by
\begin{eqnarray}
r=\tan\alpha e_1\wedge e_2, \ \ \ \ \ \ \theta=\frac{\partial}{\partial\alpha}\otimes e_3.
\end{eqnarray}
Actually, $r:\eta^*\rightarrow {\rm su}(2)\otimes {\rm su}(2)$ is a quasi-triangular dynamical r-matrix in ordinary sense where $\eta={\rm Span}\{e_3\}$ is an abelian subalgebra of ${\rm su}(2)$.

\subsection{Gauge transformations of generalized classical dynamical r-matrices}
Let $G$ be a Lie group and $r$ be a dynamical r-matrix coupled with $(U,\pi_U)$ via $\theta$ with respect to $\Omega$. Let $g:U\rightarrow G$ be a smooth map. We define a function $\langle g^*\Theta,\theta\rangle : U\rightarrow\frak g\otimes\frak g$, where $\Theta=g^{-1}dg$ is the Cartan one form on $G$ and $\langle\cdot,\cdot\rangle$ is the pairing between one forms and vector fields components of $g^*{\Theta}\in\Omega^1(U,\frak g)$ and $\theta\in\Gamma(TU\wedge\frak g)$ respectively. Similarly, with the help of the pairing between vectors and forms, we can define a function $\langle \pi_U,g^*\Theta\wedge g^*\Theta\rangle : U\rightarrow \frak g\wedge\frak g$ and a section $\langle \pi_U,g^*\Theta\rangle \in \Gamma(TU\otimes\frak g)$. For $r:U\rightarrow \frak g\wedge\frak g$ and $\theta\in\Gamma(TU\otimes\frak g)$, we define
\begin{eqnarray}
&&r^{g}:=Ad_g\otimes Ad_g(r+\widehat{\langle g^*\Theta,\theta\rangle}+\langle \pi_U, g^*\Theta\wedge g^*\Theta\rangle), \\
&&\theta^g:=Ad_g(\theta+\langle \pi_U,g^*\Theta\rangle),
\end{eqnarray}
where $\widehat{\langle g^*\Theta,\theta\rangle}:U\rightarrow \frak g\wedge\frak g$ is the antisymmetrization of $\langle g^*\Theta,\theta\rangle$.
\begin{pro}\label{gauger}
$r^{g}$ is a generalized classical dynamical r-matrix coupled with $U$ via $\theta^g$ with respect to $\Omega$.
\end{pro}
\pf Following Theorem \ref{PoissonDYBE}, given the dynamical r matrix $(\pi_U,\theta,r)$, we can construct a right invariant bivector $\pi:=\pi_U+\rho^L(\hat{\theta})+\rho^L(r)$ on $U\times G$ such that $[\pi,\pi]=-2\rho(\Omega)$. We define $U'\subset U\times G$ as the graph of the map $g:U\rightarrow G$ in $U\times G$. Then $\pi$ is a right G invariant bivector fields and $U'$ is a cross-section of the right G action. Following the argument in Theorem \ref{ReduceGDYBE}, associated to $\pi$ and $U'$ there exists a dynamical r matrix $(\pi_{U'},\theta_{U'},r_{U'})$ with respect to $\Omega$. To write down the explicit expressions of $\theta_{U'}$ and $r_{U'}$, we just need to compute the decomposition of $\pi|_{U'}$ with respect to the $TU'$ and $\rho^R(\frak g)|_{U'}$.

We introduce an isomorphism $F:U\rightarrow U'$ by $F(x)=(x,g(x))\in U'$ for any $x\in U$. A straightforward calculation gives that $$F_*\pi_U=\pi_{U'}, \ F^*\theta^{\sharp}_{U'}={\theta^g}^{\sharp} \ and \ r_{U'}\circ F=r^g.$$ It means that after identifying $U'$ with $U$ by $F$, the triple $(\pi_{U'},\theta_{U'},r_{U'})$ becomes $(\pi_U,\theta^g,r^g)$. Thus $r^g$ is naturally a classical dynamical r matrix coupled with $(U,\pi_U)$ via $\theta^g$ with respect to $\Omega$. It finishes the proof.
\qed

\begin{defi}
We call the transformation $r\rightarrow r^{g}$, $\theta\rightarrow \theta^g$ the gauge transformation for $(\pi_U,\theta,r)$ with respect to $g:U\rightarrow G$.
\end{defi}
The geometric meaning of gauge transformations of generalized dynamical $r$-matrices is illuminated in the proof of Proposition \ref{gauger}. Another interpretation is given in the following:
Let $(M,\pi_M)$ be a (quasi-)Poisson $G$-manifold and $(\pi_U,\theta,r)$ be a dynamical r-matrix with respect to $\Omega=0$(the Cartan $3$-tensor).
Given a gauge transformation $g:U\rightarrow G$, we define a transformation from $U\times M$ to itself by
\begin{eqnarray}
g\cdot(x,p)=(x,g(x)\cdot p), \ \forall p\in G.
\end{eqnarray}
\begin{pro}\label{gaugetransformation}
Let $(U\times M,\pi_r)$ and $(U\times M, \pi_{r^g})$ be the Poisson manifolds given in Theorem \ref{ReducedPoisson} associated to $(\pi_U,\theta,r)$ and $(\pi_U,\theta^g,r^g)$ respectively. Then we have
$$
\{F\circ g,G\circ g\}_r=\{F,G\}_{r^g}\circ g,
$$
for any $F,G\in C^{\infty}(U\times G)$.
\end{pro}

\subsection{Generalized classical dynamical r-matrices and Poisson groupoids}
In this subsection, we discuss the geometric interpretation of the generalized CDYB equation. Recall that in \cite{EtingofVarchenko}, Etingof and Varchenko found a geometric interpretation of the CDYB equation that
generalizes Drinfeld's interpretation of the CYB equation in terms of Poisson-Lie
groups. Namely, they constructed a so called dynamical Poisson-Lie groupoid
structure on the direct product manifold ${\eta^*}\times G\times \eta^*$,
where $\eta$ is a Lie subalgebra of $\frak g$. The CYBE can be viewed as the
special case of the generalized CYBE(see Example \ref{classicaldrmatrix}). An observation here is that $\eta^*\times G\times \eta^*$ is the Lie groupoid integrating the Lie
algebroid $T\eta^*\oplus \frak g$. Furthermore, the Poisson structure on $\eta^*\times  G\times \eta^*$ induces a Lie bialgebroid structure on $T\eta^*\oplus \frak g$. Similarly, in the case of the generalized dynamical r-matrix, we have the following theorem. Let $M$ be a manifold, $\frak g$ be a Lie algebra and $(TM\oplus\frak g,[\cdot,\cdot]_L)$ be a Lie algebroid with the anchor map given by the projection to $TM$, the bracket given by
\begin{eqnarray}
[X+A,Y+B]_L=[X,Y]+L_XB-L_YA+[A,B]_{\frak g},
\end{eqnarray}
for all $X,Y\in\Gamma(TM)$ and $A,B\in C^{\infty}(M\times\frak g)$.
\begin{thm}[Theorem 4.5, \cite{LiuXu}]\label{differential}
A solution $r$ of the generalized DYBE coupled with $(M,\pi)$ via $\theta$ induces a coboundary Lie bialgebroid structure $(TM\oplus \frak g,d_*)$ where the differential $d_*:\Gamma(\wedge^{\bullet}(TM\oplus \frak g))\rightarrow \Gamma(\wedge^{\bullet+1}  (TM\oplus \frak g))$ corresponding to the Lie algebroid structure on $T^*M\oplus \frak g^*$ is of the form
$$
d_*=[\pi_M+\hat{\theta}+r,\cdot]_L,
$$
where $[\cdot,\cdot]_L$ is the Schouten bracket on $\wedge^{\bullet}(TM\oplus\frak g)$.
\end{thm}

A solution of the generalized Poisson-Lie DYBE coupled with $(M,\pi_M)$ via $\theta$ also give a Lie bialgebroid $(TM\oplus \frak g,d_*)$ where the differential $d_*=\delta+[\pi_M+\hat{\theta}+r,\cdot]$. According to the theory of integration of Lie bialgebroids in \cite{MX}, we have
\begin{cor}\label{Poissongroupoid}
Associated to a generalized classical dynamical r-matrix coupled with $(M,\pi)$ via $\theta$, there is a Poisson Lie groupoid structure on $\mathcal{G}=M\times G\times M$ whose tangent Lie bialgebroid is $(TM\oplus\frak g,d_*)$.
\end{cor}
Thus the Poisson Lie groupoid $M\times G\times M$ gives a geometric interpretation of the generalized DYBE that generalizes Drinfeld's interpretation of the CYBE in terms of Poisson-Lie groups.

A smooth manifold $M$ is called $\mathcal{G}$-space for a Lie groupoid $(\mathcal{G}\rightrightarrows P, s,t)$ if there are two smooth maps, the moment map and the action map, $J:M\rightarrow P$ and
\begin{eqnarray*}
\alpha:\mathcal{G}\times_JM=\{(x,m)\in \mathcal{G}\times M~|~t(x)=J(m)\}\rightarrow M
\end{eqnarray*}
such that, writing $\alpha(x,m)=x\cdot m$, for all compatible $x,y\in\mathcal{G}$ and $m\in M$,

\begin{itemize}
  \item[(i)] $J(x\cdot m)=s(x);$

\item[(ii)] $(x\cdot y)\cdot m=x\cdot (y\cdot m);$

\item[(iii)] $J(m)\cdot m=m.$
\end{itemize}

Now suppose that $M$ is a $\mathcal{G}$-space. The action of $\mathcal{G}$ on $M$ is a Poisson action if its graph $\{(x,m,x\cdot m)~|~t(x)=J(m)\}$ is a coisotropic submanifold of $\mathcal{G}\times M\times \overline{M}$ \cite{Alan}. Then $M$ is called a Poisson $\mathcal{G}$-space.

From Theorem \ref{adjointPoisson} and Corollary \ref{Poissongroupoid}, associated to a dynamical $r$-matrix coupled with M via $\theta$ with respect to the Cartan 3-tensor and a quasi-Poisson G-space $N$, we have a Poisson Lie groupoid $\mathcal{G}=M\times G\times M$ and a Poisson manifold $(M\times N,\pi)$. Further, there is a natural action of the groupoid $\mathcal{G}\rightrightarrows M$ on $M\times N$ defined by $$(x,g,y)\cdot (y, p)=(x,g\cdot p)$$ for all $x,y\in M$, $g\in G$ and $p\in N$. This is a Lie groupoid action with respect to the moment map $J:M\times N\rightarrow M$ given by the natural projection. An observation here is that this action is a Poisson action.
\begin{thm}\label{Poissonaction}
$(M\times N,\pi)$ is a Poisson $\mathcal{G}$-space.
\end{thm}
To prove this theorem, we need the following results.
\begin{lem}[Theorem 3.3, \cite{liu}]\label{Liebimap}
Let $\mathcal{G}$ be a Poisson groupoid with its tangent Lie bialgebroid $(A,A^*)$. Then a Poisson manifold $(M,\pi)$ is a Poisson $\mathcal{G}-$space if and only if the vector bundle morphism from $T^*M$ to $A^*$, the dual of the infinitesimal action map, is a Lie algebroid morphism.
\end{lem}
\begin{lem}[Theorem 3.1, \cite{liu}]\label{Liemap}
Let $A^*$ be a Lie algebroid over P, and $(M,\pi)$  a Poisson manifold. Then, for the cotangent Lie algebroid $T^*M$ induced by the Poisson structure, a vector bundle morphism $\Phi:T^*M\rightarrow A^*$ over $J:M\rightarrow P$ is a Lie algebroid morphism if and only if the following two conditions hold:
\begin{itemize}
 \item[(i)] $H_{J^*f}=-\Phi^*(d_*f), \ \forall f\in C^{\infty}(P)$;

\item[(ii)] $L_{\Phi^*(X)}\pi=-\Phi^*(d_*S), \ \forall S\in\Gamma(A)$,
\end{itemize}
where $H_{J^*f}$ denotes the Hamiltonian vector field on $M$, which is defined by $H_gh=\pi(dg,dh)$ for all $g,h\in C^{\infty}(M)$. The differential $d_*$ comes from the Lie algebroid structure on $A^*$.
\end{lem}
{\bf Proof of Theorem \ref{Poissonaction}} In our case, the Poisson manifold is $M\times N$ and the Lie bialgebroid is $(TM\otimes \frak g,d_*)$. So by Lemma \ref{Liebimap} and Lemma \ref{Liemap}, we just need to prove $H_{J^*f}=-F(d_*f), \ \forall f\in C^{\infty}(M)$ and $L_{F(S)}\pi_M=-F(d_*S), \ \forall S\in\Gamma(TM\otimes \frak g)$, where the bundle map $F:TM\oplus\frak g\rightarrow T(M\times N)$ is the infinitesimal action of $M\times G\times M$ on $M\times N$, explicitly given by $F(X+e)=X+\rho(e)$ for $X\in\Gamma(TM)$ and $e\in \frak g$.

$(1)$ Following the expression of the Poisson tensor $\pi$ on $M\times N$ given in Theorem \ref{adjointPoisson}, we have for all $f\in C^\infty(M)$, $$H_{J^*f}=\pi_M^*(df)+\rho({\theta^*(df)}),$$ where $J:M\times N\rightarrow M$ is the natural projection map. On the other hand, by the definition of the differential $d_*$, $$d_*f=[\pi_M,f]+[\hat{\theta},f]+[r,f].$$ Note that $[r,f]=0$, $F([\pi_M,f])=[\pi_M,f]$, and $F([\hat{\theta},f])=F(-\theta^*(df))=-\rho({\theta^*(df)})$. Therefore, $$H_{J^*f}=-[\pi_M+\rho(\hat{\theta})+\rho(r),f]=-F(d_*f).$$

$(2)$ We assume $S=X+e\in\Gamma(TM\oplus\frak g)$. Then $$d_*S=[\pi_M+\hat{\theta}+r,X+e]=[\pi_M+\hat{\theta}+r,X]+[\pi_M+\hat{\theta}+r,e].$$ By the definition of $F:\wedge^*(TM\oplus\frak g)\rightarrow T(M\times N)$, we have that the map $F$ and the Schouten-bracket $[\cdot,\cdot]_L$ on $\wedge^*(TM\otimes\frak g)$ are commutative. As a result, $$F([\pi_M+\hat{\theta}+r,X])=[F(\pi_M)+F(\hat{\theta})+F(r),F(X)]=[\pi,X].$$ Similarly, $$F([\pi_M+\hat{\theta}+r,e])=[F(\pi_M)+F(\hat{\theta})+F(r),F(e)]=[\pi,\rho(e)].$$ Eventually, we get $-F(d_*S)=L_{F(S)}\pi$. This finishes the proof. \qed

\section{Generalized dynamical r-matrices and moduli spaces of flat connections on surfaces}
\subsection{Poisson and quasi-Poisson structures on the moduli spaces of flat connections on surfaces}
In \cite{AtiyahBott}, Atiyah-Bott introduced canonical symplectic structures on the moduli spaces of flat connections over oriented surfaces.
Let $\Sigma_{g,n}$ be a oriented surface of genus $g$ with n punctures. A convenient description of the moduli space is given by the character variety, i.e., group homomorphisms $h:\pi_1(\Sigma_{g,n})\rightarrow G$ that map the homotopy equivalence class of a loop around the $i$-th puncture to the associated conjugacy class. Two such group homomorphisms describe gauge-equivalent connections if and only if they are related by conjugation with an element of $G$. This implies that the moduli space of flat $G$-connections on $\Sigma_{g,n}$ is given by
$$
X_G(\Sigma_{g,n})=Hom_{\mathcal{C}_1,...,\mathcal{C}_n}(\pi_1(\Sigma_{g,n},G))/G=\{h\in Hom(\pi_1(\Sigma_{g,n},G)|h(m_i)\in\mathcal{C}_i\}/G.
$$
By characterising the group homomorphisms in terms of the images of the generators of $\pi_1(S_{g,n})$, the moduli space of flat connections is the set
$$
\{(M_1,...,M_n,A_1,B_1,...,A_g,B_g)\in G^{n+2g}|M_i\in \mathcal{C}_i, [B_g, A_g]\cdot\cdot\cdot [B_1, A_1]\cdot M_n\cdot\cdot\cdot M_1=1\}/G,
$$
where the quotient is taken with respect to the diagonal action of $G$ on $G^{n+2g}$. This moduli space carries a canonical Poisson structure which is obtained via Poisson reduction from the canonical symplectic structure on the space of connections on $\Sigma_{n,g}$. We use a convenient and explicit description of this Poisson structure, which is given by Poisson reduction of a Poisson structure on an enlarged ambient space $G^{n+2g}$. In this description, the $i$-th component of $G^{n+2g}$ corresponds to the holonomy along a generator of the fundamental group $\pi_1(\Sigma_{g,n})$. Thus the space of representations is embedded into $G^{n+2g}$. There are two natural operators $\nabla_R$, $\nabla_L\in \Gamma(TM\otimes{\frak g}^*)$ given for all $A\in\frak g$, $p\in G$ by $$\langle\nabla_R,A\rangle f(p):= \frac{d}{dt}|_{t=0}f(p\cdot exp(-tA)),$$ $$\langle\nabla_L,A\rangle f(p):= \frac{d}{dt}|_{t=0}f(exp(tA)\cdot p).$$ With them, we can define $2(n+2g)$ covariant differential operators in the following way:
\begin{eqnarray}
\nabla_{2i-1}=\nabla_R^{M_i}, && \nabla_{2i}=\nabla_L^{M_i} \ for \ i=1,...,n;\nonumber\\
\nabla_{n+4i-3}=\nabla_R^{A_i}, && \nabla_{n+4i-1}=\nabla_L^{A_i} \ for \ i=1,...,g;\\
\nabla_{n+4i-2}=\nabla_R^{B_i}, && \nabla_{n+4i-1}=\nabla_L^{B_i} \ for \ i=1,...,g.\nonumber
\end{eqnarray}
\begin{defi}
Let $G$ be a Lie group with Lie algebra $\frak g$. For any $r\in \frak g\otimes\frak g$, the corresponding Fock and Rosly's bivector $B_r^{n,g}\in \Gamma(\wedge^2(TG^{n+2g}))$ is defined by
\begin{eqnarray}\label{FockRosly}
B_r^{n,g}(df,dh):=\frac{1}{2}\displaystyle{\sum_{i}}\langle r,\nabla_if\wedge\nabla_i h\rangle+\sum_{i<j}\langle r,\nabla_i f\wedge\nabla_jh\rangle.
\end{eqnarray}
\end{defi}

\begin{thm}\rm{\cite{FockRosly}}
Let $\frak g$ be a Lie algebra with a non-degenerate Ad-invariant symmetric bilinear form $K(\cdot,\cdot)$. If $r\in \frak g\otimes\frak g$ is a solution of the classical Yang-Baxter equation
$$
[[r,r]]=[r_{12},r_{13}]+[r_{12},r_{23}]+[r_{13},r_{23}]=0,
$$
then $B_r^{n,g}$ defines a Poisson structure on $G^{n+2g}$.
Furthermore, when the symmetric part of $r$ coincides with the Casimir element $\kappa$, this Poisson structure induces the canonical symplectic structure on the moduli space of flat $G$-connections on $\Sigma_{g,n}$.
\end{thm}
From the expression of $B_r^{n,g}$, we see that the Poisson bracket of two functions on $G^{n+2g}$ depends only on the symmetric component of $r$ if one of the two functions is invariant under the diagonal action of $G$ on $G^{n+2g}$. So we can use the symmetric part of r and reduction procedure to describe the Poisson structure on the quotient space $G^{n+2g}/G$. In this case instead of the Poisson structure, we have a quasi-Poisson structure on $G^{n+2g}$. Actually, a main result in quasi-Poisson theory is that the standard Poisson structure on $X_G(\Sigma_{n,g})$ is given by a reduction of a canonical quasi-Poisson tensor on $M_G(\Sigma_{n,g})=\{h\in Hom(\pi_1(S_{g,n},G)|h(m_i)\in \mathcal{C}_i\}$.
\begin{thm}\label{quasiPoisson}{\rm\cite{Anton}}
Consider the quasi-Poisson manifold
$$
P_{g,n}=\mathcal{C}_1\ast...\ast \mathcal{C}_n\ast D(G)\ast...\ast D(G),
$$
where $\mathcal{C}_1$,...,$\mathcal{C}_n$ are conjugacy classes in G. Then the reductions of $P_{g,n}$ are isomorphic to the moduli spaces of flat connections on $\Sigma_{g,n}$ with the Atiyah-Bott symplectic form.
\end{thm}
A direct computation shows that the quasi-Poisson tensor on $P_{g,n}$ coincides with the restriction of $B_{\kappa}^{g,n}$ to the symplectic leaf associated to the set of conjugacy classes $\mathcal{C}_1$,...,$\mathcal{C}_n$, here $B_{\kappa}^{g,n}$ is the bivector field given in \eqref{FockRosly} with respect to the Casimir element $\kappa\in S^2\frak g$.

\subsection{GCDYB equations and moduli spaces of flat connections on surfaces}
In this subsection, we will combine the discussion in previous sections and give our main result which describes the canonical Poisson structure on the moduli spaces of flat connections on surfaces by using generalized dynamical r-matrices.

Following Theorem \ref{quasiPoisson}, the Poisson structure on $X_G(\Sigma_{g,n})$ is given by the reduction of the quasi-Poisson structure $B_{\kappa}^{g,n}$ on $P_{g,n}=\mathcal{C}_1\ast...\ast \mathcal{C}_n\ast D(G)\ast...\ast D(G)$ with respect to the simultaneous conjugation action of $G$.
\\
\\
{\bf Reduction with respect to two punctures}

We assume that there are more than 2 punctures on the surface $\Sigma_{g,n}$. We choose a cross-section of the $G$ action on $\Sigma_{g,n}$. Then the reduced Poisson structure on it is viewed to be a local model of the Poisson structure on $X_G(\Sigma_{g,n})$. We do the reduction in a ``minimal'' way, i.e., imposing gauge fixing conditions on the first two punctures in the following way. We see $P_{g,n}$ as the fusion product of $(\mathcal{C}_1\times\mathcal{C}_2,\pi_{G^2})$ and $(P_{g,n-2},B_{\kappa}^{g,n-2})$, i.e., $P_{g,n}=(\mathcal{C}_1\ast \mathcal{C}_2)\ast P_{g,n-2}$ in which $P_{g,n-2}:=\mathcal{C}_3\ast...\ast \mathcal{C}_n\ast D(G)\ast ...\ast D(G)$. Let $U$ be any cross-section of the diagonal action of G on $\mathcal{C}_1\times \mathcal{C}_2$ and $(\pi_U,\theta,r)$ the moduli space generalized dynamical r-matrix associated to $U$. By Theorem \ref{adjointPoisson}, we obtain the reduced Poisson structure on $U\times P_{g,n-2}$:
$$
\pi_{{\rm red}}=\pi_U+\rho(\hat{\theta})+\rho(r)+B_{\kappa}^{g,n-2},
$$
where $\rho:\frak g\rightarrow P_{g,n-2}$ is the infinitesimal action generated by simultaneous conjugation G action. Moreover a simple comparison gives that $\rho(r)=B_r^{g,n-2}$ as bivector fields on $U\times P_{g,n-2}$, where $B_r^{g,n-2}$ is given by \eqref{FockRosly} for $r:U\rightarrow \frak g\wedge\frak g$. As a result, $\rho(r)+B_{\kappa}^{g,n-2}=B_{r+\kappa}^{g,n-2}$, where $r$ and $\kappa$ can be seen as antisymmetric and symmetric parts of the function $r+k\in C^{\infty}(U,\frak g\otimes\frak g)$. Eventually, we get
\begin{thm}\label{reducedPoisson}
The quasi-Poisson structure $B_{\kappa}^{g,n}$ on $P_{g,n}$ induces a Poisson bracket on $U\times \mathcal{C}_3...\times\mathcal{C}_n\times G^{2g}$, which is isomorphic to the Atiyah-Bott symplectic structure and takes the following form:

$(1)$ For $f,g\in C^{\infty}(\mathcal{C}_3...\times\mathcal{C}_n\times G^{2g})$.
\begin{eqnarray}\label{Gn-2Poisson}
\{f,g\}=B_{r+\kappa}^{n-2,g}(df,dg)
\end{eqnarray}

$(2)$ For $f\in C^{\infty}(\mathcal{C}_3...\times\mathcal{C}_n\times G^{2g})$ and $\phi,\varphi\in C^{\infty}(U)$:
\begin{eqnarray}\label{thetaoperator}
&&\{f,\phi\}=\rho(\hat{\theta})(df,d\phi)\label{differentialoperator}\\
&&\{\phi,\varphi\}=\pi_U(d\phi,d\varphi),\label{UPoisson}
\end{eqnarray}
\end{thm}
Note that the original Fock and Rosly's bivector $B_r^{n,g}$ on $G^{n+2g}$ is related to a classical r-matrix which is a constant element in $\frak g\otimes\frak g$. Here we describe the Poisson tensor on the quotient space by a dynamical r-matrix which is a map from $U$ to $\frak g\wedge \frak g$.
\begin{cor}
Let $\mathcal{G}$ be the Poisson Lie groupoid associated to the moduli space dynamical r matrix $(\pi_U,r,\theta)$, then $(U\times \mathcal{C}_3...\times\mathcal{C}_n\times G^{2g},\{\cdot,\cdot\})$ carries a Poisson $\mathcal{G}$ action.
\end{cor}
{\bf Gauge fixing and classical dynamical $r$-matrices in $\ISO(2, 1)$-Chern-Simons theory.}

In \cite{Meusburger}, Meusburger and Sch$\rm\ddot{o}$nfeld obtained classical dynamical $r$-matrices by considering gauge fixing in $\ISO(2,1)$-Chern-Simons theory. Now, we interpret these classical dynamical $r$-matrices as moduli
space dynamical $r$-matrices corresponding to $G=\ISO(2,1)$.

First, let us give the required notations.
We denote by $e_0=(1,0,0)$, $e_1 = (0,1,0)$, $e_2 = (0,0,1)$ the standard basis of $\mathbb R^3$. All indices run from $0$
to $2$ and are raised and lowered with the three-dimensional Minkowski metric $\eta =diag(1,-1,-1)$. By $\varepsilon_{abc}$ we denote the totally antisymmetric tensor
in three dimensions with the convention $\varepsilon_{012}=1$.

The Poincar$\rm\acute{e}$ group in 3-D is the semidirect product $\ISO(2, 1) =\SO_+(2, 1)\ltimes\mathbb R^3$ of the proper orthochronous Lorentz group $\SO_+(2,1)$ and the translation group $\mathbb R^3$.
We parameterize elements of $\ISO(2, 1)$ as
\begin{eqnarray*}
(u,\mathbf{a})=(u,0)\cdot(1,-\mathbf{j})=(u,-Ad(u)\mathbf{j}) \ with \ u\in SO_+(2,1),~\mathbf{j,a}\in\mathbb R^3.
\end{eqnarray*}
The corresponding coordinate functions $\{j^a\}_{a=0,1,2}$ is given by $$j^a : \ISO(2, 1)\rightarrow \mathbb R,\quad (u, -Ad(u)\mathbf{q})\rightarrow q^a.$$
We fix a set of generators $\{J_a\}_{a=0,1,2}$ of $\frak {so}(2,1)$ such that the Lie bracket takes the form $[J_a,J_b]=\varepsilon_{ab}^{ \ \ c}J_c$.  Then a basis of the Lie algebra $\frak {iso}(2,1)$ is given by the basis $\{J_a\}_{a=0,1,2}$ of $\frak {so}(2,1)$ together with a
basis $\{P_a\}_{a=0,1,2}$ of the abelian Lie algebra $\mathbb R^3$.

The moduli space of flat $G$-connections can be viewed as a constrained system in the sense of Dirac \cite{Dirac}. From the expression for the moduli space of flat connections, it is apparent that the moduli space is obtained from $P_{g,n}$ by imposing a group-valued constraint that arises from the defining relation of the fundamental group $\pi_1(\Sigma_{g,n})$.

In the case of $G=\ISO(2,1)$, the group-valued constraint from the defining relation of the fundamental group $\pi_1(\Sigma_{g,n})$ can be viewed as a set of constraints in the Dirac gauge fixing formalism for the Fock-Rosly bracket. The associated gauge transformations which they generate via the Poisson bracket are given by the diagonal action of $\ISO(2,1)$ on $\ISO(2,1)^{n+2g}$.

A choice of gauge fixing conditions for the constraints \eqref{constraint} is investigated in \cite{Meusburger}. 
These gauge fixing conditions implement the quotient by $\ISO(2, 1)$ and restrict the variables $M_1, M_2$ in such a way that for all points $(M_1, . . . , B_g)\in \Sigma = C^{-1}(0)$, the components $M_1, M_2 \in \ISO(2, 1)$ are determined uniquely by two real parameters $\psi$ and $\alpha$ given in terms of the Lorentzian and translational components of the product $M_2\cdot M_1=(u_{12}, -Ad(u_{12})\mathbf{j}_{12})$ as
\begin{eqnarray}
\psi=f(Tr(u_{12})), \ \ \ \ \ \alpha=g(Tr(u_{12} ))Tr(j^a_{12} J_a \cdot u_{12} )+h(Tr(u_{12} )),
\end{eqnarray}
where $f, g \in C^{\infty}(\mathbb R)$ are arbitrary diffeomorphisms and $h \in C^{\infty}(R)$. This allows us to identify the constraint surface $\Sigma = C^{-1}(0)$ with a subset of $\mathbb R^2 \times \ISO(2, 1)^{n-2+2g}$, where the $\mathbb R^2$ is parameterized by $(\psi,\alpha)$ and $\ISO(2,1)^{n-2+2g}$ by $(M_3,...,B_g)$.

We choose a set of constraint functions of the form \eqref{constraint}, \eqref{Constraint1} and \eqref{Constraint2} that satisfy the requirements of Dirac gauge fixing. Following Theorem 4.5 in \cite{Meusburger}, there exist maps $\mathbf{q}_{\psi}, \mathbf{q}_{\alpha}, \mathbf{q}_{\delta},m:\mathbb R^2\rightarrow\mathbb R^3$, $V:\mathbb R\rightarrow {\rm Mat}(3,\mathbb R)$ such that the associated Dirac bracket is given in terms of them. On the other hand, the Dirac gauge fixing corresponds to the constraint functions is equivalent to choose a cross-section of the $\ISO(2,1)$ action on $\mathcal{C}_1\times\mathcal{C}_2$ which is the locus of the six functions $\mathcal{C}_{i=1,...,6}$. Thus by Theorem \ref{reducedPoisson}, there is a moduli space dynamical $r$-matrix corresponding to this cross-section. We give the following theorem which rewrites the main result in \cite{Meusburger} using our language.
\begin{thm}
The moduli space dynamical classical r-matrix $(\pi,\theta,r)$ corresponding to the Dirac gauge fixing procedure is given by
\begin{eqnarray}
&&\pi=0, \ \ \ \ \ \theta=q_{\alpha}^a\frac{\partial}{\partial \alpha}\otimes J_a+q_{\psi}^a\frac{\partial}{\partial\psi}\otimes P_a+q_{\delta}^a\frac{\partial}{\partial\alpha}\otimes P_a,\label{ISOtheta}\\
&&r=-V^{bc}(\psi)(P_b \otimes J^c-J^c \otimes P_b)+\varepsilon^{bcd}m_d(\psi,\alpha)P_b\otimes P_c\label{ISOr}.
\end{eqnarray}
Moreover, the induced Poisson bracket takes the following form:\\
for any $f,g\in C^{\infty}(\ISO(2,1)^{n-2+2g})$,
\begin{eqnarray}\label{Gn-2Poisson}
\{\alpha,\psi\}=0, \ \ \{\alpha,f\}=\rho(\hat{\theta})(d\alpha,df), \ \
\{f,g\}=B_{r+\kappa}^{n-2,g}(df,dg)
\end{eqnarray}
where $\kappa=P_a\otimes J^a$.
\end{thm}

Given a map $p:\mathbb{R}^2\rightarrow \ISO(2,1)$, we consider smooth maps
\begin{eqnarray*}
\Phi^p :\mathbb R^2 \times \ISO(2, 1)^{n-2+2g}\rightarrow \mathbb R^2\times \ISO(2, 1)^{n-2+2g},\\
(\psi,\alpha,M_3,...,B_g)\rightarrow (\psi,\alpha,Ad_pM_3,...,Ad_pB_g).
\end{eqnarray*}
By Proposition \ref{gaugetransformation}, we have
\begin{cor}\label{gaugetransformation}
Let $\{\cdot,\cdot\}_D$ be the bracket given in \eqref{Gn-2Poisson} with respect to $(\theta,r)$. Then for all $F,G\in C^{\infty}(\mathbb R^2\times \ISO(2,1)^{n-2+2g})$, we have
\begin{eqnarray}
\{F\circ \Phi^p,G\circ\Phi^p\}_D=\{F,G\}_D^p\circ\Phi^p,
\end{eqnarray}
where $\{\cdot,\cdot\}_D^p$ is the bracket given in \eqref{Gn-2Poisson} with respect to $(\theta^p,r^p)$, the gauge transformation of $(\theta,r)$ via $p:\mathbb R^2\rightarrow \ISO(2,1)$.
\end{cor}
Particularly, the map $p=(g,-Ad(g)\mathbf{t})$ satisfying $\partial_{\alpha}g=\partial_{\alpha}^2\mathbf{t}=0$ is called dynamical Poincar$\rm\acute{e}$ transformation in \cite{Meusburger}. Actually, dynamical $r$-matrices from different gauge fixing conditions subject to extra conditions $(a)$ and $(b)$ are related by dynamical Poincar$\acute{e}$ transformations. In this special case, Corollary \ref{gaugetransformation} turns to be Lemma 5.1 in \cite{Meusburger}.

A simplified standard set of dynamical $r$-matrices from the Dirac gauge fixing in $\ISO(2,1)$-Chern-Simons theory is given explicitly in \cite{Meusburger}. This set of solutions corresponds to special gauge fixing condition which is motivated by its direct physical interpretation in the application to the Chern-Simons formulation of $(2+1)$-gravity.


\begin{thebibliography}{}
%\cite{Alexandrov:1995kv}


\bibitem{Anton}
A. Alekseev, Y. Kosmann-Schwarzbach and E. Meinrenken, {\em Quasi-Poisson manifolds}, Canadian J. Math. 54 (2002), 3-29.

\bibitem{AMM}
A. Alekseev, A. Malkin, and E. Meinrenken, {\em Lie group valued moment maps}, J. Differential Geom, 48 (1998) 445–495.

\bibitem{AtiyahBott}
M. Atiyah and R. Bott, {\em The Yang-Mills equations over Riemann surfaces}, Philosophical Transactions of the Royal Society of London, Series A, Mathematical and Physical Sciences, 54 (1983) 523-615.

\bibitem{Dirac}
P. A. M. Dirac, {\em Generalized Hamiltonian dynamics}, Canadian Journal of Mathematics, 2 (1950) 129–148.

\bibitem{EnriquezEtingof}
B. Enriquez, P. Etingof and I. Marshall, {\em Quantization of some Poisson-Lie dynamical r-matrices and Poisson homogeneous spaces},
Quantum groups, 135–175, Contemp. Math., 433, Amer. Math. Soc., Providence, RI, 2007.

\bibitem{EtingofVarchenko}
J. Donin and A. Mudrov, {\em Dynamical Yang-Baxter equation and quantum vector bundles}, Comm. Math. Phys. 254 (2005), no. 3, 719-760.

\bibitem{AB}
G. Felder, {\em Conformal field theory and integrable systems associated to elliptic curves}, Proc. ICM Z$\ddot{u}$rich, (1994), 1247-1255.

\bibitem{FockRosly}
V. V. Fock and A. A. Rosly, {\em Poisson structure on moduli of flat connections on Riemann surfaces and the r-matrix}, Moscow Seminar in Mathematical Physics, Amer. Math. Soc. Transl. Ser. 2, 191 (1999), 67–86.

\bibitem{Gold}
W. Goldman, {\em The symplectic nature of fundamental groups of surfaces}, Adv. in Math., 54(2):200–225, 1984.

\bibitem{liu}
L.-G. He, Z.-J. Liu and D.-S. Zhong, {\em Poisson actions and Lie bialgebroid morphisms}, Quantization, Poisson brackets and beyond (Manchester, 2001), 235–244, Contemp. Math., 315, Amer. Math. Soc., Providence, RI, 2002.

\bibitem{LiuXu}
Z.-J. Liu and P. Xu, {\em The local structure of Lie bialgebroids}, Lett. Math. Phys. 61 (2002), no. 1, 15–28.

\bibitem{Meusburger}
C. Meusburger and T. Sch$\rm\ddot{o}$nfeld, {\em Gauge fixing and classical dynamical r-matrices in ISO(2,1)-Chern-Simons theory}, Comm. Math. Phys. 327 (2014), no. 2, 443–479.

\bibitem{MX}
K.C.H. Mackenzie and P. Xu, {\em Lie bialgebroids and Poisson groupoids}, Duke Math. J. 73 (1994), 415–452.

\bibitem{Alan}
A. Weinstein, {\em Coisotropic calculus and Poisson groupoids}, J. Math. Soc. Japan. 40 (1988), 705–727.

\bibitem{PingXu}
P. Xu, {\em Triangular dynamical r-matrices and quantization}, Adv. Math. 166 (2002), no. 1, 1–49.

\end{thebibliography}
\end{document}